\documentclass[12pt,a4paper]{article}

\usepackage{graphicx}
\usepackage{amssymb}
\usepackage{amsmath}
\usepackage{caption}
\usepackage{subcaption}
\usepackage[margin=0.8in]{geometry} 
\usepackage{subfloat}
\usepackage{color}

\title{Numerical solution of scattering problems using a Riemann--Hilbert formulation}

\author{
Stefan G. Llewellyn Smith$^{1,2}$ and Elena Luca$^{1}$}

\date{}

\newcommand{\ci}{\mathrm{i}}
\newcommand{\Cc}{\mathcal{C}}

\newcommand{\dd}{\mathrm{d}}
\newcommand{\ee}{\mathrm{e}}

\usepackage{MnSymbol}
\usepackage{bigints}

\begin{document}

\maketitle

\begin{center}
$^{1}$Department of Mechanical and Aerospace Engineering \\
Jacobs School of Engineering, UCSD \\
La Jolla, CA 92093-0411, USA.
\end{center}

\begin{center}
$^{2}$Scripps Institution of Oceanography, UCSD \\
La Jolla, CA 92039-0213, USA.

\vskip 0.1truein 
{\tt sgls@ucsd.edu} \\
{\tt elouca@eng.ucsd.edu}
\end{center}

\vskip 0.5truein

\begin{center}
{\bf Abstract}
\end{center}
\noindent
A fast and accurate numerical method for the solution of scalar and matrix Wiener--Hopf problems is presented. The Wiener--Hopf problems are formulated as Riemann--Hilbert problems on the real line, and a numerical approach developed for these problems is used. It is shown that the known far-field behaviour of the solutions can be exploited to construct numerical schemes providing spectrally accurate results. A number of scalar and matrix Wiener--Hopf problems that generalize the classical Sommerfeld problem of diffraction of plane waves by a semi-infinite plane are solved using the approach.

\vfill\eject

\section{Introduction}

The Wiener--Hopf (WH) method was originally developed to solve integral equations and mixed boundary value problems \cite{WienerHopf:1931,LawrieAbrahams:2007}. Standard references include \cite{Noble:1988,Weinstein:1969,DanieleZich:2014}. In its direct formulation introduced by Jones \cite{Jones:1952}, the WH method, in combination with the Fourier transform, reduces the solution of a boundary-value problem (in e.g.~electrodynamics, acoustics, elasticity, etc\ldots) to the
problem of solving a functional equation by finding functions analytic in the upper and lower complex half-planes by factorization. The typical WH functional equation (with $\alpha$ the transform variable) is
\begin{equation}
A(\alpha) \Phi_{+}(\alpha) + B(\alpha) \Phi_{-}(\alpha) + C(\alpha)=0, \label{WH problem}
\end{equation}
valid in the strip $a_-<\text{Im}(\alpha)<a_+$, where the functions $\Phi_{+}(\alpha)$, $\Phi_{-}(\alpha)$ are, respectively, analytic in $a_-<\text{Im}(\alpha)$ and $\text{Im}(\alpha)<a_+$. Further, $A(\alpha)$,
$B(\alpha)$ and $C(\alpha)$ are known functions of $\alpha$. In addition, the behaviour of $\Phi_{+}(\alpha)$, $\Phi_{-}(\alpha)$ for large $|\alpha|$ is given by the behaviour of physical variables near the origin.

While it is straightforward to apply the WH method to a scalar problem for which an explicit solution formula based on Cauchy's integral formula exists \cite{Noble:1988}, there is no general method for carrying out the decomposition for matrix functions. A survey of constructive methods for factorization problems is given by Rogosin \& Mishuris \cite{RogosinMishuris:2016}. The decomposition can be carried out for matrices of special form, following the ideas of Khrapkov \cite{Khrapkov:1971}, Hurd \cite{Hurd:1976}, Daniele \cite{Daniele:1978} and Rawlins \& Williams \cite{RawlinsWilliams:1981}, but for other matrices, the only general semi-analytical approach that has been put forward is the Pad\'e decomposition method of Abrahams \cite{Abrahams:1997}, in which the matrices are approximated by rational functions, permitting a decomposition. The method seems to work well in applications \cite{ADLS1,ADLS2} but a lot of manipulation of the equations and algebra is required. A recent study by Kisil \cite{Kisil:2018} proposes an iterative method for triangular matrix problems with exponential factors.

WH problems are related to Riemann--Hilbert (RH) problems. The former connect boundary values of sectionally analytic functions in a common strip of analyticity, while the latter couple these on a contour. The RH problem asks for the construction of a function $\Phi(\alpha)$ that is analytic everywhere in the complex plane except along a given oriented contour $\Gamma$ where it has a prescribed jump
\begin{equation}
A(\alpha) \Phi_{+}(\alpha)+B(\alpha) \Phi_{-}(\alpha)+C(\alpha)=0, \qquad \alpha \in \Gamma, \label{RH problem}
\end{equation}
where $\Phi_{+}(\alpha)$ and $\Phi_{-}(\alpha)$ are the representations of $\Phi(\alpha)$ on the $+$ and $-$ sides of the contour, respectively and $A(\alpha)$, $B(\alpha)$ and $C(\alpha)$ are known functions on $\Gamma$. The functional form \eqref{RH problem} is identical to \eqref{WH problem}, but it is now valid along a given oriented contour rather than a strip. Many RH problems arise in the context of singular integral equations, but RH problems have also been connected to random matrix theory, nonlinear special functions, nonlinear wave equations, and other problems. Olver \cite{Olver:2012} and Trogdon \& Olver \cite{TO:2015} have recently developed accurate and efficient numerical algorithms for the solution of RH problems.

The aim of this study is to present fast and accurate numerical schemes for the solution of scalar and matrix WH problems by exploiting the links between the WH and RH problems. The idea is to solve the associated RH problems, adapting the methods of Olver \cite{Olver:2012} and Trogdon \& Olver \cite{TO:2015} to take into account the known far-field behaviour of the solutions. In particular, the approach adopted is to use rational mappings with multiple inverses, extending the M\"obius mappings discussed in \cite{Olver:2012,TO:2015}, to account for the $\mathcal{O}(\alpha^{-1/2})$ decay of the solutions for large $|\alpha|$. We focus here on problems of diffraction of plane waves by a semi-infinite plane that produce scalar and matrix WH problems. For the acoustic scattering problems considered, the conditions on $A(\alpha), B(\alpha), C(\alpha)$ guaranteeing existence and uniqueness of solutions are satisfied \cite{Gakhov:1966,TO:2015}. These schemes are also applicable to other more complicated diffraction problems, such as the diffraction of plane waves by two identical strips \cite{Shanin:2003}.

In \S\,2 we present the numerical scheme for the solution of Riemann--Hilbert problems using rational mappings. In \S\,3 we give error measures that we use to verify our numerical solutions. This is followed by implementation of the numerical scheme for a number of scalar and matrix problems for diffraction by a semi-infinite plane (\S\,4). In \S\,5 we compute the far-field directivity, $D(\theta)$, which is defined via
\begin{equation}
\phi(r,\theta) \sim D(\theta) \frac{\mathrm{e}^{{\rm i}kr}}{\sqrt{r}}, \quad \text{as} \quad r\rightarrow \infty, \label{directivity}
\end{equation}
where $\phi$ is the diffracted field, $k$ is the acoustic wave number and $(r,\theta)$ are polar coordinates. Finally we conclude and discuss further applications in \S\,6.

\section{Numerical solution of Riemann--Hilbert problems}

Following Olver \cite{Olver:2012} and Trogdon \& Olver \cite{TO:2015}, we present the numerical scheme for the solution of RH problems \eqref{RH problem}. Although the method can be used to solve RH problems over any contour $\Gamma$ whose individual pieces can be mapped from the unit interval $\mathbb{I}=[-1,1]$, here we take $\Gamma=\mathbb{R}$, since this is relevant to WH problems.

The aim is to find a function $\Phi(\alpha)$ analytic everywhere in the complex plane except along $\mathbb{R}$ where it has a prescribed jump of the form \eqref{RH problem}. We represent the function $\Phi(\alpha)$ by
\begin{equation}
\Phi(\alpha)=\sum_{k=0}^{n-1} \check U_k T_k(x)=\sum_{k=0}^{n-1} \check U_k T_k(M^{-1}(\alpha)) \label{Phiexpansion}
\end{equation}
for $\alpha \in \mathbb{R}$ and by analytic continuation off
$\mathbb{R}$. Here $T_k(x)$ is the usual $k$th Chebyshev polynomial and $M: \mathbb{I}=[-1,1] \mapsto \mathbb{R}$ is a mapping
\begin{equation}
\alpha=M(x), \quad x=M^{-1}(\alpha).
\end{equation}
In this study, we shall consider rational mappings $M$ with multiple inverses, extending the M\"obius mappings discussed in \cite{Olver:2012,TO:2015}. (The notation $M^{-1}$ refers to mapping back to $\mathbb{I}$; subscripts for $M^{-1}$ will indicate other branches.) The function $\Phi(\alpha)$ can be scalar, vector or matrix, and the coefficients $\check U_k$ are of the same nature. The coefficients $\check U_k$ can be found from the function values via
\begin{equation}
\check U_k = \sum_{q=1}^{n} F_{kq} \Phi(\alpha_q), \label{coefficients}
\end{equation}
where $F_{kq}$ are the elements of the Chebyshev operator $\mathcal{F}$. We seek a collocation method in which the function $\Phi$ and the relation \eqref{RH problem} are evaluated at points $\alpha_{q}=M(x_{q})$ along the contour $\Gamma=\mathbb{R}$, where $\{x_{q}|q=1,\dots,n\}$ are the Chebyshev points of the second kind:
\begin{equation}
{\bf x}^{\mathbb{I}}=(x_{1},\dots,x_{n})=\left(-1,\cos\frac{(n-2)\pi}{n-1},\dots,\cos\frac{\pi}{n-1},1\right). \label{collocationpts}
\end{equation}

\subsection{Collocation method}

The Cauchy operator is defined by
\begin{equation}
\mathcal{C} \Phi(\alpha)=\frac{1}{2\pi {\rm i}} \int_{\Gamma} {\frac{\Phi(\alpha')}{\alpha'-\alpha} \mathrm{d}\alpha'}. \label{Cauchyop}
\end{equation}
The Cauchy operators $\mathcal{C}^{\pm}$ are defined as the limiting values of the Cauchy operator \eqref{Cauchyop} as $\alpha$ approaches the oriented contour $\Gamma=\mathbb{R}$ from the $+$ and $-$ sides, i.e.~$\mathcal{C}^{+} \Phi=\Phi_{+}$ and $\mathcal{C}^{-}\Phi=\Phi_{-}$. To construct a collocation method for solving \eqref{RH problem}, we must compute the Cauchy matrices $C^{\pm}$ at points $\alpha_{q}=M(x_{q})$, $x_{q} \in {\bf x}^{\mathbb{I}}$, along $\Gamma$. (The notation $\mathcal{C}^{\pm}$ refer to the operators; $C^{\pm}$ are matrices.) Since $C^{+}-C^{-}=I$, it is sufficient to construct $C^{+}$.

Consider
\begin{equation}
G(\alpha) = \sum_{k=0}^{n-1} \check U_k [{\cal C}T_k](M^{-1}(\alpha))= \sum_{k=0}^{n-1} \check U_k S_k(M^{-1}(\alpha)),
\end{equation}
where $S_k \equiv {\cal C}T_k$. This is the Cauchy transform of $\Phi(\alpha)$, ${\cal C}\Phi(\alpha)$, up to a constant which requires analysis for large $|\alpha|$. Next, we evaluate $G(\alpha)$ at the point $\alpha_{p}=M(x_p)$, where $x_{p}$ is a Chebyshev point in \eqref{collocationpts}, giving
\begin{equation}
G(\alpha_p) = \sum_{k=0}^{n-1} \check U_k S_k(M^{-1}(\alpha_p)).
\end{equation}
Substituting in for the coefficients $\check U_k$ \eqref{coefficients}, we can write this as a linear operator acting on the function values:
\begin{equation}\label{G_ap}
G(\alpha_p) = \sum_{k=0}^{n-1} \sum_{q=1}^{n} F_{kq} S_k(M^{-1}(\alpha_p)) \Phi(\alpha_q) = \sum_{q=1}^{n} C_{pq} \Phi(\alpha_q),
\end{equation}
where we have defined
\begin{equation}\label{Cpqinitial}
C_{pq} = \sum_{k=0}^{n-1} F_{kq} [{\cal C}T_k](M^{-1}(\alpha_p)) =  \sum_{k=0}^{n-1} F_{kq} S_k(M^{-1}(\alpha_p)) .
\end{equation}
In the next section, we consider rational mappings with $d$ inverses and, therefore, we must write
\begin{equation}\label{Cpqmultiple}
C_{pq} = \sum_{j=1}^{d} C^{(j)}_{pq},
\end{equation}
where $C^{(j)}_{pq}$ is given by \eqref{Cpqinitial} with $M^{-1}$ replaced by the $j$-th inverse $M_j^{-1}$. When the contour $\Gamma$ is unbounded with the endpoints $\pm 1$ of ${\bf x}^{\mathbb{I}}$ the preimages of $|\alpha|=\infty$, the above form of ${\cal C}\Phi(\alpha)$ needs to take into account the limiting behaviour for large $|\alpha|$ \cite{Olver:2012}; this will be discussed in the next section.

\subsection{Rational mappings and Cauchy matrices}

Olver \cite{Olver:2012} and Trogdon \& Olver \cite{TO:2015} considered mappings $M$ that are M\"{o}bius transformations and also discussed the extension to polynomials of degree $d$ with $d$ inverses. In what follows, we extend their formulation to rational functions $M$. This allows us to represent functions that have $\mathcal{O}(\alpha^{-1/2})$ decay for large $|\alpha|$. The use of the Chebyshev expansion \eqref{Phiexpansion} and subsequent matrix operators imply spectral convergence of the numerical scheme \cite{Trefethen:2013}.

\subsubsection{2-to-1 mapping}

\begin{figure}
\centering
\includegraphics[scale=1]{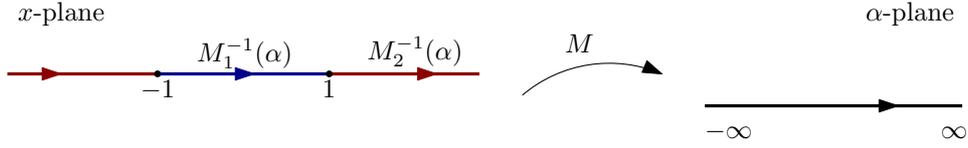}
\caption{The rational 2-to-1 mapping given by \eqref{map1}. The inverses $M^{-1}_{j}(\alpha)$, $j=1,2$ are given by \eqref{inverses2to1}.}
\label{figmap2to1}
\end{figure}

We consider the rational mapping
\begin{equation}\label{map1}
\alpha = M(x)=\frac{x}{1-x^2}
\end{equation}
which has two inverses:
\begin{equation}
M^{-1}_{1,2}(\alpha) = \frac{-1 \pm \sqrt{1+4\alpha^2}}{2\alpha} \quad \Rightarrow \quad M^{-1}_{1}(\alpha)=x, \quad M^{-1}_{2}(\alpha)=-x^{-1}. \label{inverses2to1}
\end{equation}
The points $x=\pm 1$ are mapped to $\pm \infty$ and, therefore, the images of $M_{j}^{-1}(\alpha)$, $j=1,2$ are both $\mathbb{R}$. A schematic is given in figure \ref{figmap2to1}. The rational mapping \eqref{map1}, which is an odd function of $x$, has the desirable property that it maps the unit $x$-interval to the entire real axis, as needed for WH problems, with $\infty$ having two preimages ($x=\pm 1$).

For each inverse, we construct the Cauchy transform matrix. Following \cite{Olver:2012, TO:2015}, the contribution of $M^{-1}_{1}(\alpha)=\mathbb{I}$ to $C^{+}$ is given by
\begin{equation}
C^{+}_{1}=\frac{1}{{\rm i}\pi} \left( \begin{array}{ccc} \frac{\log2}{2}+\frac{{\rm i}\pi}{2} & & \\ & -2\text{ diag}(\text{arctanh } T_{\downarrow}^{-1}({\bf x}^{\mathbb{I}})) & \\ & & -\frac{\log2}{2}+\frac{{\rm i}\pi}{2} \end{array} \right)+\frac{1}{{\rm i}\pi} F^{-1}P F, \label{C1plus}
\end{equation}
where $F$ is the transform matrix from the function values at ${\bf x}^{\mathbb{I}}$ to the coefficients, $P$ is an almost Toeplitz matrix and
\begin{equation}
T_{\downarrow}^{-1}(x)=x-{\rm i}\sqrt{1-x}\sqrt{1+x}.
\end{equation}
Explicit expressions for $F$ and $P$ are given in \cite{Olver:2012, TO:2015}.

The contour $M^{-1}_{2}(\alpha)=\mathbb{R}\setminus \mathbb{I}$ is unbounded and connected to $\mathbb{I}$ at the points $\pm 1$. The contribution of this contour to $C^{+}$ is
\begin{equation}
C^{+}_{2}=\left(\begin{array}{c} \mu^{R} \\ \Psi^{z}_{n} \\ \mu^{L} \end{array}\right) F, \label{2to1C2}
\end{equation}
where $\mu^{L/R}$ are row vectors (where $L$ is identified with $-1$ and $R$ with $+1$) and $\Psi^{z}_{n}$ is an $n \times n$ matrix with the  first and last rows removed, so that the matrix in \eqref{2to1C2} is of dimension $n \times n$. Both are defined in \cite{Olver:2012, TO:2015} and are associated with the Cauchy transforms of the Chebyshev basis evaluated at points
\begin{equation}
z=T^{-1}_{+}(M^{-1}_{2}(\alpha)), \quad \text{where} \quad T^{-1}_{+}(x)=x-\sqrt{x-1}\sqrt{1+x}.
\end{equation}
As in \cite{Olver:2012, TO:2015}, the orientation of the unbounded contour $M^{-1}_{2}(\alpha)$ (connected to $\mathbb{I}$ at $\pm 1$) implies that the first and last rows of the matrix in \eqref{2to1C2} must be the row vectors $\mu^{R}$ and $\mu^{L}$, respectively.

The final form of the Cauchy matrix $C^{+}$ for the rational mapping \eqref{map1} is given by
\begin{equation}
C^{+}=C^{+}_{1}+C^{+}_{2}- {\bf 1}_{n\times n} ~\text{diag}(\mu^{L}+\mu^{R}) ~ F. \label{Cplus2to1}
\end{equation}
The final term is required to ensure that $C^+$ has the correct limiting behaviour \cite{Olver:2012}. Our numerical scheme takes into account the behaviour of the transform functions $\Phi_\pm$ for large $|\alpha|$. Their decay at infinity implies that the first and last rows in the matrices $C^{\pm}$ are irrelevant in our calculations.

However, this mapping does not possess the desired structure at the endpoints, since, e.g.~for $x\sim 1$,
\begin{equation}
\alpha \sim \frac{1}{2(1-x)}
\end{equation}
and similarly for $x \sim -1$. For the diffraction problems we consider, the functions $\Phi_{\pm}$ are of the form $a_1 \alpha^{-1/2} + a_2 \alpha^{-1} + a_3 \alpha^{-3/2} + \cdots$ and, therefore, their representations in this case will not be spectrally convergent Chebyshev series in $x$.

\subsubsection{4-to-1 mapping}

\begin{figure}
\centering
\includegraphics[scale=1]{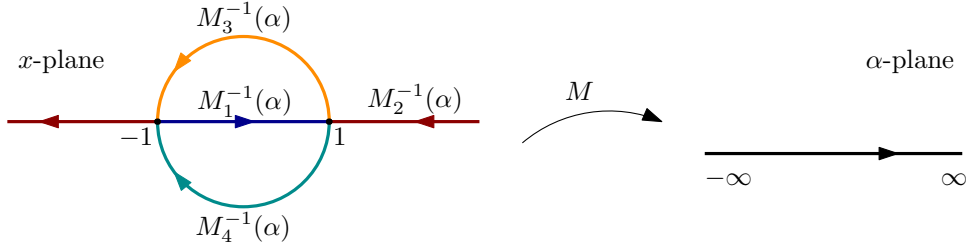}
\caption{The rational 4--to--1 mapping given by \eqref{map2}. The inverses $M^{-1}_{j}(\alpha)$, $j=1,2,3,4$ are given by \eqref{inverses4to1a}--\eqref{inverses4to1b}.}
\label{figmap4to1}
\end{figure}

We introduce the rational mapping
\begin{equation}\label{map2}
\alpha = M(x)=\frac{x+x^3}{(1-x^2)^2}
\end{equation}
which has four inverses:
\begin{equation}
M^{-1}_{1}(\alpha)=x, \quad M^{-1}_{2}(\alpha)=x^{-1}, \label{inverses4to1a}
\end{equation}
and the solutions of a quadratic given by
\begin{equation}
\quad M^{-1}_{3,4}(\alpha)=\frac{-c(x)\pm\sqrt{c(x)^2-4}}{2}, \quad \text{where} \quad c(x)=x+x^{-1}-\alpha^{-1}. \label{inverses4to1b}
\end{equation}
A schematic is given in figure \ref{figmap4to1}. Again, the points $x=\pm 1$ are mapped to $\pm \infty$, but now the mapping will also be able to capture exactly the far-field behaviour of the functions $\Phi_{\pm}$. The critical point is that, for $x \sim 1$,
\begin{equation}
\alpha \sim \frac{1}{2(1-x)^2}
\end{equation}
and similarly for $x \sim -1$. One can hence represent functions of the form $a_1 \alpha^{-1/2} + a_2 \alpha^{-1} + a_3 \alpha^{-3/2} + \cdots$ as spectrally convergent Chebyshev series in $x$.

Again, for each one of the inverses, we construct the Cauchy matrix. The contour $M^{-1}_{1}(\alpha)$ is identical to that presented above for the 2-to-1 mapping and, therefore, its contribution to $C^{+}$, $C^{+}_{1}$, is again given by \eqref{C1plus}.

The contour $M^{-1}_{2}(\alpha)=\mathbb{R}\setminus \mathbb{I}$ is unbounded and connected to $\mathbb{I}$ at the points $\pm 1$. The orientation of the contour is shown in figure \ref{figmap4to1}. Its contribution to $C^{+}$ is given by \eqref{2to1C2} evaluated at $ z=T^{-1}_{+}(M^{-1}_{2}(\alpha))$. Note that, consistent with \cite{Olver:2012, TO:2015}, the orientation of the unbounded contour $M^{-1}_{2}(\alpha)$ (connected to $\mathbb{I}$ at $\pm 1$) implies that the first and last rows of the matrix must be the row vectors $\mu^{L}$ and $\mu^{R}$, respectively.

The contours $M^{-1}_{3}(\alpha)=\{\mathrm{e}^{{\rm i}t}| t \in [0,\pi]\}$, $M^{-1}_{4}(\alpha)=\{\mathrm{e}^{{\rm i}t}| t \in [-\pi,0]\}$ are bounded and connected to $\mathbb{I}$ at the points $\pm 1$ (their orientation is shown in figure \ref{figmap4to1}). The contributions of these contours to $C^{+}$ are given again by \eqref{2to1C2}, evaluated at $z=T^{-1}_{+}(M^{-1}_{3}(\alpha))$ and $z=T^{-1}_{+}(M^{-1}_{4}(\alpha))$, respectively. Again, consistent with \cite{Olver:2012, TO:2015}, the orientation of the bounded contours $M^{-1}_{3}(\alpha)$, $M^{-1}_{4}(\alpha)$ (connected to $\mathbb{I}$ at $\pm 1$) implies that the first and last rows of the matrix must be the row vectors $\mu^{R}$ and $\mu^{L}$, respectively.

The final form of the Cauchy matrix $C^{+}$ for the rational mapping \eqref{map2} is given by
\begin{equation}
C^{+}=C^{+}_{1}+C^{+}_{2}+C^{+}_{3}+C^{+}_{4}- {\bf 1}_{n\times n} ~\text{diag}~2(\mu^{L}+\mu^{R}) ~ F. \label{Cplus4to1}
\end{equation}
The fifth term has been subtracted to ensure that $C^+$ has the correct limiting behaviour (Olver \cite{Olver:2012}). Our numerical scheme takes into account the behaviour of the transform functions $\Phi_\pm$ for large $|\alpha|$. Their decay at infinity implies that the first and last rows in the matrices $C^{\pm}$ are again irrelevant for our calculations.

\subsection{Rotation, scaling and evaluation \label{scalingrotation}} 

Once the matrix $C^{+}$ has been constructed (using either the 2-to-1 or 4-to-1 mapping), $C^{-}$ follows from the relation $C^{+}-C^{-}=I$. Therefore, we have a collocation method for the RH equation \eqref{RH problem} which can be written as
\begin{equation}
[A(\alpha_p) C^{+}+B(\alpha_p) C^{-}]\Phi(\alpha_p)+C(\alpha_p)=0
\label{RHnote}
\end{equation}
for $\alpha_{p}=M(x_p)$, $p=1,\dots, n$, and $x_{p}$ given in \eqref{collocationpts}. The WH problem requires the contour to lie in the strip of analyticity which encloses $\Gamma=\mathbb{R}$. For numerical purposes it is helpful to take the contour $\Gamma$ far from singularities, for example by taking a rotation of the real axis by angle $-\chi$, e.g.~evaluating \eqref{RHnote} along
\begin{equation}\label{rotatedalpha}
\alpha_p^{\text{rotated}}=\alpha_p \mathrm{e}^{-{\rm i}\chi},
\end{equation}
provided that no branch points are traversed. A RH problem can be then formulated along $\mathbb{R} \mathrm{e}^{-{\rm i}\chi}$; this poses no problem given the decay properties of the functions $\Phi_{+}$ and $\Phi_{-}$ and the fact that the rotation matrix providing \eqref{rotatedalpha} is invertible. Specifically, if we multiply \eqref{RHnote} by a rotation matrix, then in order for \eqref{RHnote} still to hold, we require functions $A(\alpha)$ and $B(\alpha)$ to be well-behaved in the region between the initial and rotated contour.

We also note that for the diffraction problems to be considered in the present study, there is one associated length scale, the diffraction parameter $k$. This implies that the mapping $M$ should be rescaled incorporating the length scale $k$. This is not pursued here, since we shall be considering the case $k=1$. For more complicated problems, e.g.~the problem of Wickham and Abrahams \cite{Wickham:1995,Abrahams:1997} discussed below, there are 2 associated length scales ($1$ and $k$), which should be taken into account in the scaling of the mapping.

To evaluate $\mathcal{C} \Phi$ at a point $\alpha=M(x)$ off the contour $\Gamma$, we use \eqref{Phiexpansion} and obtain $\mathcal{C} \Phi(\alpha)=\sum_{k=0}^{n-1} \check U_k \mathcal{C} T_k(x)$, where the coefficients $\check U_k$ are known e.g.~from the numerical solution of the Riemann--Hilbert problem. To compute $\mathcal{C} T_k(x)$ for each of the inverses $M^{-1}_{j}(\alpha)$, we compute a row vector $\Psi^z$, analogous to the matrix $\Psi_n^z$, corresponding to the Cauchy transform of the Chebyshev polynomials at a single point $\alpha=M(x)$.

\section{Error estimate}

In the following sections we apply the numerical scheme presented above to solve various scalar and matrix WH problems. To validate our results, we compare computed values for the sectionally analytic functions $\Phi_{\pm}$ against those given by exact solutions. We use the error estimate
\begin{equation}
E^{r}_{n} = \left[\int_{-1}^{1} {|Q(x)-Q_n(x)|^{r} \mathrm{d}x} \right]^{1/r} \label{error}
\end{equation}
($r > 1$; for $r = \infty$ the estimate becomes the maximum error), where $Q(x)$ is the exact solution and $Q_n(x)$ is the numerical value at the collocation points. Note that \eqref{error} is a measure of the error in the mapped unit interval $\mathbb{I}$ in the $x$-plane. This integral can be evaluated using Clenshaw--Curtis quadrature \cite{ClenshawCurtis:1960} using the values of the integrand at the Chebyshev points. We write
\begin{equation}
E^{r}_{n} \approx \left[ \sum_{q=1}^{n} w_{q} |Q(x_{q})-Q_n(x_{q})|^{r} \right]^{1/r},
\end{equation}
where $\{x_{q}|q=1,\dots,n\}$ are the Chebyshev points \eqref{collocationpts} and $w_{q}$ are weights that can be found in e.g.~\cite{Trefethen:2013, DavisRabinowitz:1984}.

Alternatively, the error estimate in the variable $\alpha$ can be used:
\begin{equation}
\mathcal{E}^{r}_{n} = \left|\int_{-\infty}^{\infty} {|R(\alpha)-R_n(\alpha)|^{r} \mathrm{d}\alpha} \right|^{1/r}, \label{erroralpha}
\end{equation}
where $R(\alpha)$ is the exact solution and $R_n(\alpha)$ is the numerical value for $n$ collocation points in the $\alpha$-plane. The error function $\mathcal{E}^{r}_{n}$ can be equivalently written as
\begin{equation}
\mathcal{E}^{r}_{n} = \left| \int_{-1}^{1} {|R(\alpha(x))-R_n(\alpha(x))|^{r} ~\frac{\mathrm{d}\alpha}{\mathrm{d} x}~ \mathrm{d}x} \right|^{1/r}.
\end{equation}
This can also be computed using Clenshaw--Curtis quadrature.

\section{Diffraction by a semi-infinite plane}

\begin{figure}
\centering
\includegraphics[scale=1]{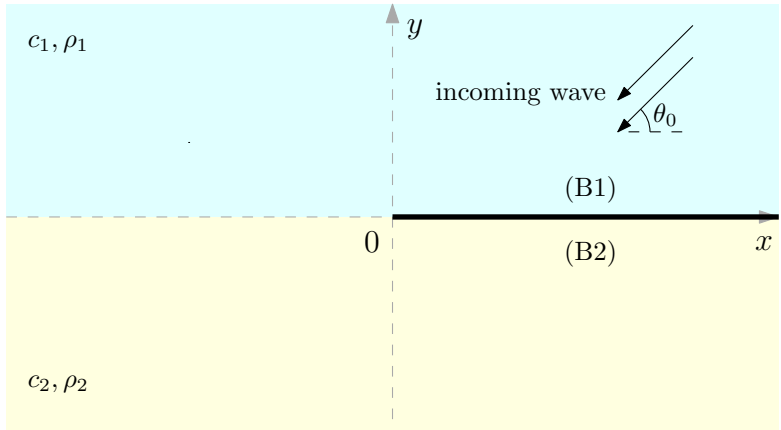}
\caption{Diffraction of plane wave by a semi-infinite plane.}
\label{fig1}
\end{figure}

We analyse the problem of diffraction of a plane wave by a half-plane, where the plane is taken to lie along the positive real axis $x>0$, $y=0$. A schematic is shown in figure \ref{fig1}. The parameters $\rho_{1}, \rho_{2}$ and $c_{1}, c_{2}$ are the density and sound speeds in the upper (medium 1) and lower (medium 2) half-planes respectively. The boundary conditions along the top and bottom sides of the semi-infinite plane along the positive real axis are denoted by (B1) and (B2).

The governing equations for the fields $\phi_{1}$ and $\phi_{2}$ in the upper and lower half-planes, respectively, are given by
\begin{equation}
\label{Helmholtz12}
\nabla^{2} \phi_{1}+k_{1}^{2} \phi_{1}=0, \qquad \nabla^{2} \phi_{2}+k_{2}^{2} \phi_{2}=0,
\end{equation}
where $k_1=\omega/c_1$, $k_2=\omega/c_2$ with $\omega$ the frequency. We decompose the total field $\phi_t$ into an incident wave $\phi_{\text{inc}}$ and a scattered field,
with $\phi_{\text{inc}}$ given by
\begin{equation}\label{phiinc}
\phi_{\text{inc}}=\text{exp}(-{\rm i}kx \cos\theta_0-{\rm i}ky\sin\theta_0), 
\end{equation}
where $k$ is the acoustic wave number and $\theta_0$ is the angle of the incident wave.
Note that there is also an underlying time dependence, $\mathrm{e}^{-{\rm i} \omega t}$, that is only needed to determine branch cuts later.

The general form of the boundary conditions (B1) and (B2) on the semi-infinite plane $x>0$, $y=0^{\pm}$ is
\begin{equation}
\lambda \phi_{t} + \frac{\partial \phi_{t}}{\partial y}=0, \label{impedancebcgeneral}
\end{equation}
for $\lambda \in \mathbb{C}$ which can be different on the two sides of the semi-infinite plane. An acoustically hard plane corresponds to $\lambda = 0$ while an acoustically soft plane corresponds to $|\lambda|\rightarrow \infty$. 

Near the origin, we write 
\begin{equation}
\phi_{t}(r,\theta) \sim r^{a} h(\theta), \label{localpolar}
\end{equation}
using polar coordinates $x = r \cos \theta$, $y = r \sin \theta$, where $a$ and $h(\theta)$ are determined from the boundary and interface conditions.

Table 1 shows a list of diffraction problems by a
half-plane with boundary conditions (B1) and (B2) along the top and
bottom sides of the semi-infinite plane respectively, ratio of sound
speeds $c_{1}/c_{2}$, ratio of densities $\rho_{1}/\rho_{2}$, the
value of $a$ determining the local behaviour near the origin
\eqref{localpolar}, as well as references to the original papers
treating each case. Different notations and conventions are used in
these references, and we follow them for ease of comparison, rather
than recasting all problems in a single unified notation. Since our goal is to show the effectiveness of the numerical method presented in \S\,2, we do not pursue an exhaustive parameter study. We note that the numerical scheme is, in general, robust with regards to parameter changes, except for the cases $\theta_0 \rightarrow \pi/2$ and impedance parameter $|\lambda| \rightarrow 0$ in \eqref{impedancebcgeneral} where convergence becomes slower. In the former case, this reduces the scope for using the angle of rotation used to separate the contour from the branch points. In the latter case, the behaviour near the origin which governs convergence properties could possibly change from the singularity structure considered in our numerical scheme near the origin ($a=1/2$).

\begin{table}
\centering
\begin{tabular}{|c|c|c|c|c|c|}
\hline
(B1) & (B2) & $c_{1}/c_{2}$ & $\rho_{1}/\rho_{2}$ & $a$ & References\\
\hline\hline
hard & hard & 1 & 1 & 1/2 & Sommerfeld \cite{Noble:1988, Sommerfeld1896} \\
soft & soft & 1 & 1 & 1/2 & Sommerfeld \cite{Noble:1988, Sommerfeld1896} \\
soft & hard & 1 & 1 & 1/4 & Rawlins \cite{Rawlins:1975} \\
impedance$_1$ & impedance$_1$ & 1 & 1 & 1/2 & Senior \cite{Noble:1988, Senior:1952} \\
impedance$_1$ & impedance$_2$ & 1 & 1 & 1/2 & Hurd \cite{Hurd:1976}, Barton \& Rawlins \cite{BartonRawlins:1999} \\
soft & hard & $\neq 1$ & $\neq 1$ & $\pi^{-1} \tan^{-1}{\mu}$ & Wickham \cite{Wickham:1995}, Abrahams \cite{Abrahams:1997} \\
\hline
\end{tabular}
 \label{tablecases}
\caption{Diffraction by a semi-infinite plane and boundary conditions on top and bottom sides of the semi-infinite boundary $x>0$, $y=0$. Note that in some cases \cite{Abrahams:1997, Noble:1988, Rawlins:1975} the semi-infinite plane lies along the negative real axis. The parameter $\mu$ is defined to be $\mu=\sqrt{\rho_2/\rho_1}$.}
\end{table}

We define full and half-range Fourier transforms in $x$ according to
\begin{equation}
\Phi(\alpha,y) = \int_{-\infty}^\infty \phi(x,y) \ee^{\ci\alpha x}
\,\dd x=\Phi_{+}(\alpha,y)+\Phi_{-}(\alpha,y), \label{Fouriertransform}
\end{equation}
with
\begin{equation}
\Phi_+(\alpha,y) = \int_0^\infty \phi(x,y) \mathrm{e}^{\ci\alpha x} \,\dd x,
\qquad \Phi_-(\alpha,y) = \int_{-\infty}^0 \phi(x,y) \mathrm{e}^{\ci\alpha x}
\,\dd x.
\end{equation}
This definition of the Fourier transform will be used in the next sections for the Sommerfeld problem (a) and the Senior problem (b). Note that Noble's analysis\cite{Noble:1988}, which will be followed to formulate both problems, uses \eqref{Fouriertransform} with a different normalising factor. For the Hurd problem (c), a different definition of the Fourier transform pair is used (following Hurd \& Przezdziecki \cite{HurdPrze:1981}).

\subsection{The Sommerfeld problem \cite{Noble:1988, Sommerfeld1896,Crighton2012modern}}

First, we revisit the Sommerfeld problem. This problem can be formulated as a scalar Wiener--Hopf problem and admits an exact solution, since the associated kernel function can be factorized explicitly into a product of an upper and a lower analytic function.

For the diffraction problem originally solved by Sommerfeld \cite{Sommerfeld1896}, we have $c_{1}=c_{2}$, $\rho_{1}=\rho_{2}$, with hard-hard boundary conditions along $x<0$, $y=0^{\pm}$. (Following \cite{Noble:1988}, the semi-infinite boundary is taken to lie along the negative real axis unlike in the rest of this paper). Then
\begin{equation}
\label{BCSom}
\text{(B1)}: \quad \frac{\partial \phi_{t}}{\partial y}=0, \quad x<0, ~y=0^{+}, \qquad
\text{(B2)}: \quad \frac{\partial \phi_{t}}{\partial y}=0, \quad x<0, ~y=0^{-}.
\end{equation}
We can write $\phi_{t}=\phi_{\text{inc}}+\phi$, where $\phi_{\text{inc}}$ is an incident wave of the form \eqref{phiinc} and $\phi$ is the scattered potential everywhere in the complex plane. Here we take $-\pi/2<\theta_0<\pi/2$; for $\theta_0$ outside this range, the analysis must be adapted. Noble \cite{Noble:1988} presents the following analysis stating that $0<\theta_0<\pi$, but a careful observation of the analyticity properties used indicates that the problem should really be considered separately for the two cases $0<\theta_0<\pi/2$ and $\pi/2<\theta_0<\pi$, even though the final result for directivity is the same. On the edge of the plate, $\phi$ is bounded and $\partial\phi/\partial y = O(x^{-1/2})$. Substitution of $\phi_t$ into \eqref{BCSom}, gives the boundary conditions for the scattered field $\phi$:
\begin{equation}\label{scatteredBCSom}
\frac{\partial \phi}{\partial y} = \ci k \sin\theta_0 \exp{(-\ci k x \cos \theta_0)}, \quad x<0, ~y=0^{\pm}.
\end{equation}

\subsubsection{Exact solution \cite{Noble:1988}}
It is convenient to assign a small positive imaginary part to $k \mapsto k+{\rm i}\epsilon$, $\epsilon>0$, to improve convergence of subsequent Fourier integrals (and then let $\epsilon \rightarrow 0$).

Since $\dd\Phi/\dd y = \Phi'$ is continuous everywhere, we write
\begin{equation}\label{expansionPhiSom}
\Phi(\alpha,y)=\begin{cases}
A(\alpha) \mathrm{e}^{-\gamma y}, \quad & y \geq 0, \\
-A(\alpha) \mathrm{e}^{\gamma y}, \quad & y \leq 0,
\end{cases}
\end{equation}
where $\gamma(\alpha)= (\alpha^2 - k^2)^{1/2}$ which has branch cuts from $\pm k$ to $\infty$ in the first and third quadrants, as required by causality. The function $A(\alpha)$ has to be determined from the boundary conditions. It follows that
\begin{equation}
\label{Som3eqns}
\Phi_+(0) + \Phi_-(0^{+}) =  A(\alpha), \qquad \Phi_+(0) + \Phi_-(0^{-}) =  -A(\alpha), \qquad \Phi_+'(0) + \Phi_-'(0) =  -\gamma A(\alpha).
\end{equation}
Expressions with one argument, e.g.~$\Phi_+(0)$, will always refer to the value of, in this case, $\Phi_+(\alpha,y)$ at $y=0$. Taking the half-range Fourier transform of the boundary condition \eqref{scatteredBCSom} gives
\begin{equation}\label{Phidminus0}
\Phi_-'(0) = \frac{k \sin{\theta_0}}{\alpha - k \cos{\theta_0}}.
\end{equation}
Elimination of the function $A(\alpha)$ from \eqref{Som3eqns} and use of \eqref{Phidminus0} gives
\begin{equation}\label{RHSom}
\Phi_+(0) = - S_-, \qquad \Phi_+'(0) + \frac{k
  \sin{\theta_0}}{\alpha - k \cos{\theta_0}} =  - \gamma
D_-,
\end{equation}
where
\begin{equation}
\Phi_-(0^{+}) - \Phi_-(0^{-}) = 2 D_-, \qquad \Phi_-(0^{+}) + \Phi_-(0^{-}) = 2 S_-.
\end{equation}
The two equations in \eqref{RHSom} hold in the strip $-k_i < \text{Im } \alpha < k_i \cos{\theta}_0$
where $k_i=\text{Im } k$. Next, we divide the second equation in \eqref{RHSom} by $(\alpha+k)^{1/2}$:
\begin{equation}\label{RHSom2}
\frac{\Phi_+'(0)}{(\alpha + k)^{1/2}} + \frac{k
  \sin{\theta_0}}{(\alpha + k)^{1/2} (\alpha - k
  \cos{\theta_0})} =  - (\alpha - k)^{1/2} D_-.
\end{equation}
The first term on the left-hand side is upper analytic, while the right-hand side is lower analytic. Applying an additive splitting to the second term gives
\begin{equation}
\frac{k
  \sin{\theta_0}}{\alpha - k \cos{\theta_0}} \left(
    \frac{1}{(\alpha + k)^{1/2}} - \frac{1}{(k + k
      \cos{\theta_0})^{1/2}} \right) + \frac{k \sin{\theta_0}}{(k + k
    \cos{\theta_0})^{1/2} (\alpha - k \cos{\theta_0})} =H_{+}(\alpha)+H_{-}(\alpha).
\end{equation}
Therefore, \eqref{RHSom2} can be written in the form
\begin{equation}
J(\alpha) = (\alpha + k)^{-1/2} \Phi_+'(0) + H_+(\alpha) = - (\alpha -
k)^{1/2} D_- - H_-(\alpha).
\end{equation}
The function $J(\alpha)$ is regular in the whole plane by analytic
continuation and is bounded from edge conditions. Hence, by Liouville's
theorem, $J(\alpha) = 0$, implying
\begin{equation}
\label{Phidplus0}
\Phi'_{+}(0)=-(\alpha+k)^{1/2} H_{+}(\alpha), \qquad D_{-}=-(\alpha-k)^{-1/2} H_{-}(\alpha).
\end{equation}
The function $A(\alpha)$ follows from \eqref{Som3eqns}.

\subsubsection{Numerical solution using the RH formulation}

We now solve the hard-hard Sommerfeld problem numerically by expressing it as a RH problem on the (rotated) real line.
Recall that
\begin{equation}\label{Numsol1}
\qquad \Phi_+'(0)  + \gamma D_- =- \frac{k \sin{\theta_0}}{\alpha - k \cos{\theta_0}}.
\end{equation}
We divide by $\gamma$ so that
\begin{equation}
\gamma^{-1} \Cc^+\Phi + \Cc^- \Phi =- \frac{k \sin\theta_0}{\gamma(\alpha - k\cos\theta_0)},
\end{equation}
where $\Phi_+'(0) = \Cc^+\Phi$ and $D_- = \Cc^- \Phi$. This becomes the matrix problem $M \Phi = N$, where
\begin{equation}
\label{Numsol2}
M=\text{diag}[1/\gamma] ~ C^{+} + C^{-}, \qquad N=-\frac{k \sin\theta_0}{\gamma(\alpha-k \cos\theta_0)}.
\end{equation}
The Cauchy operator $C^{+}$ is given by \eqref{Cplus2to1} or \eqref{Cplus4to1} and $C^{-}=C^{+}-I$. The sectionally analytic functions $\Phi'_{+}(0)$ and $D_{-}$ can be computed by applying the operators $C^{+}$ and $C^{-}$ to $\Phi$.

\begin{figure}
\centering
{\hspace{-3.2em}
\includegraphics[scale=0.3]{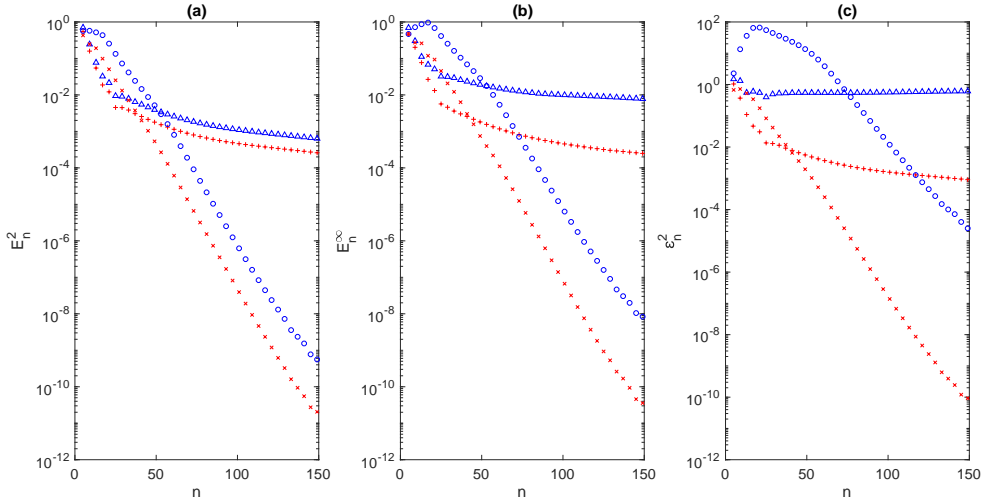}}
\caption{Error estimates (a) $E^{2}_{n}$, (b) $E^{\infty}_{n}$ and (c) $\mathcal{E}^{2}_{n}$ for the sectionally analytic functions $\Phi'_{+}(0)$ and $D_{-}$ (hard-hard Sommerfeld problem) using the 2-to-1 and 4-to-1 mappings, as functions of $n$. Parameters are $\chi=\pi/4$, $k=1$, $\theta_0=\pi/5$. Symbols: $\Phi'_{+}(0)$ using 2-to-1 mapping ({\color{blue} $\medtriangleup$}), $D_{-}$ using 2-to-1 mapping ({\color{red} $+$}), $\Phi'_{+}(0)$ using 4-to-1 mapping ({\color{blue} $\medcircle$}), $D_{-}$ using 4-to-1 mapping ({\color{red} $\times$}).}
\label{errorestimates}
\end{figure}

Figure \ref{errorestimates} shows the error estimates $E^{2}_{n}$, $E^{\infty}_{n}$, $\mathcal{E}^{2}_{n}$ for $\Phi'_{+}(0)$ and $D_{-}$ using the 2-to-1 and 4-to-1 mappings, as functions of $n$ (the number of collocation points). We observe that spectral convergence is obtained using the 4-to-1 mapping, since, as discussed earlier, this mapping is able to capture exactly the far-field behaviour of the sectionally analytic functions. On the other hand, the error estimate using the 2-to-1 mapping is much larger, e.g.~$E^{2}_n=\mathcal{O}(10^{-3})$ using $n=100$; this is to be expected, since the numerical scheme using this mapping does not incorporate the square root singularity near the origin which governs convergence properties \cite{Boyd}.

\subsection{Senior's problem \cite{Noble:1988, Senior:1952}}

We take $c_{1}=c_{2}$, $\rho_{1}=\rho_{2}$ with impedance boundary conditions along $x>0$, $y=0^{\pm}$:
\begin{equation}
\label{BCSenior}
\text{(B1)}: \quad \phi_{t}-{\rm i} S \frac{\partial \phi_{t}}{\partial y}=0, \quad x>0, ~y=0^{+}, \qquad \text{(B2)}: \quad \phi_{t}+{\rm i} S \frac{\partial \phi_{t}}{\partial y}=0, \quad x>0, ~y=0^{-},
\end{equation}
where $S$ is the impedance parameter. Again, we can write $\phi_{t}=\phi_{\text{inc}}+\phi$, where $\phi$ is the scattered potential everywhere in the complex plane and $\phi_{\text{inc}}$ is the incident wave \eqref{phiinc}, but now with $\pi/2<\theta_0<3\pi/2$ (this is omitted in \cite{Noble:1988}).

\subsubsection{Analysis \cite{Noble:1988}}
We write
\begin{equation}\label{SeniorPhi}
\Phi(\alpha,y)= \begin{cases}
A(\alpha) \mathrm{e}^{-\gamma y}, \quad & y \geq 0, \\
B(\alpha) \mathrm{e}^{\gamma y}, \quad & y \leq 0,
\end{cases}
\end{equation}
where $\gamma= (\alpha^2 - k^2)^{1/2}$ and $A(\alpha)$, $B(\alpha)$ have to be determined from the boundary conditions. This gives (using the same notation as previously):
\begin{equation}
\begin{split}\label{Senior4eqns}
\Phi_{+}(0^{+})+\Phi_{-}(0)=A(\alpha), \qquad &\Phi'_{+}(0^{+})+\Phi'_{-}(0)=-\gamma A(\alpha), \\
\Phi_{+}(0^{-})+\Phi_{-}(0)=B(\alpha), \qquad &\Phi'_{+}(0^{-})+\Phi'_{-}(0)=\gamma B(\alpha).
\end{split}
\end{equation}
Elimination of $A(\alpha)$, $B(\alpha)$, substitution of $\phi_{t}=\phi_{\text{inc}}+\phi$ in \eqref{BCSenior} and the definition of the Fourier transform yield
\begin{align}
\Phi'_{-}(0)+\gamma \Phi_{-}(0)&=-(1+{\rm i}S \gamma) \Phi'_{+}(0^{+})+ \frac{{\rm i} \gamma (1-S k \sin\theta_0)}{\alpha-k \cos\theta_0}, \label{systemWH1} \\
\Phi'_{-}(0)-\gamma \Phi_{-}(0)&=-(1+{\rm i}S \gamma) \Phi'_{+}(0^{-})- \frac{{\rm i} \gamma (1+S k \sin\theta_0)}{\alpha-k \cos\theta_0}. \label{systemWH2}
\end{align}
These equations can be expressed in matrix form as
\begin{equation}\label{Seniormatrix}
\left[ \begin{array}{cc} 1+{\rm i}S \gamma & 0 \\ 0 & 1+{\rm i}S \gamma \end{array} \right] \Phi_{+}+\left[ \begin{array}{cc} 1 & \gamma \\ 1 & -\gamma \end{array} \right] \Phi_{-}= \frac{{\rm i} \gamma}{\alpha-k \cos\theta_0}  \left[ \begin{array}{c} 1-S k \sin\theta_0 \\ -1- S k \sin\theta_0 \end{array} \right],
\end{equation}
where
\begin{equation}
\Phi_{+}=\left[ \begin{array}{c} \Phi'_{+}(0^{+}) \\ \Phi'_{+}(0^{-}) \end{array} \right], \qquad \Phi_{-}=\left[ \begin{array}{c} \Phi'_{-}(0) \\ \Phi_{-}(0) \end{array} \right].
\end{equation}
Addition and subtraction of \eqref{systemWH1} and \eqref{systemWH2} give two independent scalar WH problems:
\begin{align}
2 \Phi'_{-}(0)&=-(1+{\rm i}S \gamma) [\Phi'_{+}(0^{+})+\Phi'_{+}(0^{-})] - \frac{2{\rm i} \gamma S k \sin\theta_0}{\alpha-k\cos\theta_0}, \label{Seniorscalar1} \\
2 \gamma \Phi_{-}(0)&=-(1+{\rm i}S \gamma) [\Phi'_{+}(0^{+})-\Phi'_{+}(0^{-})] + \frac{2{\rm i} \gamma}{\alpha-k\cos\theta_0}. \label{Seniorscalar2}
\end{align} 

\subsubsection{Numerical solution using the RH formulation}

To solve this problem numerically, we express it as a RH problem. Use of \eqref{Seniormatrix} results in the matrix RH problem $M \Phi = N$, where
\begin{equation}
M=\left[ \begin{array}{cc} \text{diag}(1/\gamma+{\rm i}S) & 0 \\ 0 & \text{diag}(1/\gamma+{\rm i}S) \end{array} \right] \left[ \begin{array}{cc} C^{+} & {\bf 0} \\ {\bf 0} & C^{+} \end{array} \right] + \left[ \begin{array}{cc} \text{diag}(1/ \gamma) & I \\ \text{diag}(1/ \gamma) & -I \end{array} \right] \left[ \begin{array}{cc} C^{-} & {\bf 0} \\ {\bf 0} & C^{-} \end{array} \right]
\end{equation}
and
\begin{equation}
N= \frac{{\rm i}}{\alpha-k \cos\theta_0}  \left[ \begin{array}{c} 1-S k \sin\theta_0 \\ -1- S k \sin\theta_0 \end{array} \right].
\end{equation}
The functions $\Phi_{+}$ and $\Phi_{-}$ can be computed by applying the operators $C^{+}$ and $C^{-}$ to $\Phi$. One can also solve the two scalar WH problems given by \eqref{Seniorscalar1} and \eqref{Seniorscalar2}  numerically.

\begin{figure}
\centering
{\hspace{-3.2em}
\includegraphics[scale=0.3]{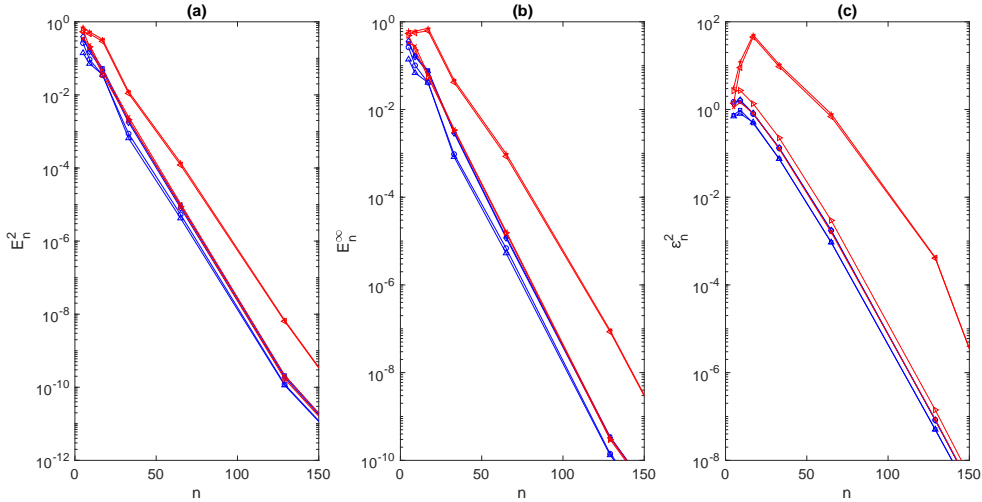}}
\caption{Error estimates for the Senior problem (a) $E^{2}_{n}$, (b) $E^{\infty}_{n}$ and (c) $\mathcal{E}^{2}_{n}$ for the sectionally analytic functions $\Phi'_{+}(0^{+}), \Phi'_{+}(0^{-}), \Phi'_-(0), \Phi_-(0)$ using the 4-to-1 mapping and both scalar and matrix formulations, as functions of $n$. Parameters are $\chi=\pi/4$, $k=1$, $\theta_0=5\pi/6$, $S=\csc(\pi/5)$. Symbols: $\Phi'_{+}(0^+)$ using the scalar form ({\color{blue} $\medtriangleup$}), $\Phi'_{+}(0^+)$ using the matrix form ({\color{blue} $\medcircle$}), $\Phi'_{+}(0^-)$ using the scalar form ({\color{blue} $\medsquare$}), $\Phi'_{+}(0^-)$ using the matrix form ({\color{blue} $\meddiamond$}), $\Phi'_{-}(0)$ using the scalar form ({\color{red} $+$}), $\Phi'_{-}(0)$ using the matrix form ({\color{red} $\medtriangleleft$}), $\Phi_{-}(0)$ using the scalar form ({\color{red} $\times$}), $\Phi_{-}(0)$ using the matrix form ({\color{red} $\medtriangleright$}).}
\label{errorestimatesSenior}
\end{figure}

Figure \ref{errorestimatesSenior} shows the error estimates $E^{2}_{n}$, $E^{\infty}_{n}$, $\mathcal{E}^{2}_{n}$ for $\Phi'_{+}(0^{+}), \Phi'_{+}(0^{-}), \Phi'_{-}(0), \Phi_{-}(0)$ using the 4-to-1 mapping and both scalar and matrix formulations, as functions of $n$. Again, spectral convergence is observed, since the 4-to-1 mapping is able to capture exactly the far-field behaviour of the sectionally analytic functions. The exact solutions used to compute the error estimates are provided from the highest resolution ($n = 257$) numerical solutions.

\subsection{Hurd's problem \cite{Hurd:1976, HurdPrze:1981, Rawlins:1984, BartonRawlins:1999}}

Hurd \cite{Hurd:1976} presented a method to deal with matrix WH problems with symmetric branch cuts at $\pm k$ and isolated singularities in the upper (or lower) half-plane. The analysis is based on the so called Wiener--Hopf--Hilbert method which involves the transformation of WH equations into a pair of coupled Hilbert equations that can be solved using Muskhelishvili's theory \cite{Muskhelishvili:1953}. The problem was revisited by Hurd \& Przezdziecki \cite{HurdPrze:1981} who carried out a more detailed analysis. Hurd's problem \cite{Hurd:1976} is similar to that analysed by Senior \cite{Noble:1988, Senior:1952}, but with different impedance boundary conditions on top and bottom sides of the semi-infinite plane:
\begin{equation}
\label{BCHurd}
\text{(B1)}: \quad \frac{\partial \phi_{t}}{\partial y}+{\rm i}k S_{1} \phi_{t}=0, \quad x>0, ~y=0^{+}, \qquad \text{(B2)}: \quad \frac{\partial \phi_{t}}{\partial y}-{\rm i}k S_{2} \phi_{t}=0, \quad x>0, ~y=0^{-}, 
\end{equation}
where $S_{j}=\sin \theta_{j}$, $0 \leq \theta_{j} \leq \pi/2$, for $j=1,2$. Again, we write $\phi_{t}=\phi_{\text{inc}}+\phi$, where $\phi_{\text{inc}}$ is the incident wave \eqref{phiinc} with $-\pi/2<\theta_0<\pi/2$ (see figures 1 and 2 in \cite{HurdPrze:1981}) and $\phi$ is the scattered potential everywhere in the complex plane. This problem was also analysed by Rawlins \cite{Rawlins:1984} using a different approach and, later, by Barton \& Rawlins \cite{BartonRawlins:1999} who modified the Wiener--Hopf--Hilbert method to consider the case when surface waves can propagate along the two sides of the semi-infinite plane.

We write
\begin{equation}\label{HurdPhi}
\Phi(\alpha,y)= \begin{cases}
A(\alpha) \mathrm{e}^{-\beta y}, \quad & y \geq 0, \\
B(\alpha) \mathrm{e}^{\beta y}, \quad & y \leq 0,
\end{cases}
\end{equation}
where $\beta(\alpha)=(k^{2}-\alpha^{2})^{1/2}$ such that $\beta(0)=k$ and $A(\alpha)$, $B(\alpha)$ have to be determined from the boundary conditions.
Following Hurd \& Przezdziecki \cite{HurdPrze:1981}, the above problem leads to a matrix WH problem of the form:
\begin{equation}\label{Hurdequations1}
\Phi_{+}-\frac{1}{2} \left[ \begin{array}{cc} \beta+k S_{1} & 1+k S_{1} /\beta \\ -\beta-k S_{2} & 1+k S_{2} /\beta \end{array} \right] \Phi_{-}=\frac{1}{2\pi {\rm i}} \frac{k}{\alpha+k \cos\theta_0}  \left[ \begin{array}{c} S_{1}-S_{0} \\ S_{2}+S_{0} \end{array} \right],
\end{equation}
where
\begin{equation}\label{Hurdequations2}
\Phi_{+}=\left[ \begin{array}{c} U_{1} \\ U_{2} \end{array} \right], \qquad \Phi_{-}=\left[ \begin{array}{c} L_{1} \\ L_{2} \end{array} \right],
\end{equation}
and $S_{0}=\sin\theta_{0}$. Note that the definition of the Fourier transform pair in \cite{HurdPrze:1981} differs from that given by \eqref{Fouriertransform}.

\subsubsection{Numerical solution using the RH formulation}

To solve this problem numerically, we express it as a RH problem. We write $M \Phi = N$, where
\begin{equation}
M=\left[ \begin{array}{cc} \text{diag}(1/\beta) & {\bf 0} \\ {\bf 0} & I \end{array} \right]  \left[ \begin{array}{cc} C^{+} & {\bf 0} \\ {\bf 0} & C^{+} \end{array} \right] - \frac{1}{2} \left[ \begin{array}{cc} \text{diag}(1+k S_{1}/\beta) & \text{diag}(1+k S_{1} /\beta) \\ \text{diag}(-1-k S_{2}/\beta) & \text{diag}(1+k S_{2}/\beta) \end{array} \right] \left[ \begin{array}{cc} C^{-} & {\bf 0} \\ {\bf 0} & C^{-} \end{array} \right]
\end{equation}
and
\begin{equation}
N=\frac{k}{2\pi{\rm i}} \left[ \begin{array}{c} (S_{1}-S_{0})/[\beta(\alpha+k \cos\theta_0)] \\ (S_{2}+S_{0})/(\alpha+k \cos\theta_0) \end{array} \right].
\end{equation}
The functions $\Phi_{+}$, $\Phi_{-}$ can be computed by applying the operators $C^{+}$ and $C^{-}$ to $\Phi$.

\begin{figure}
\centering
{\hspace{-3.2em}
\includegraphics[scale=0.3]{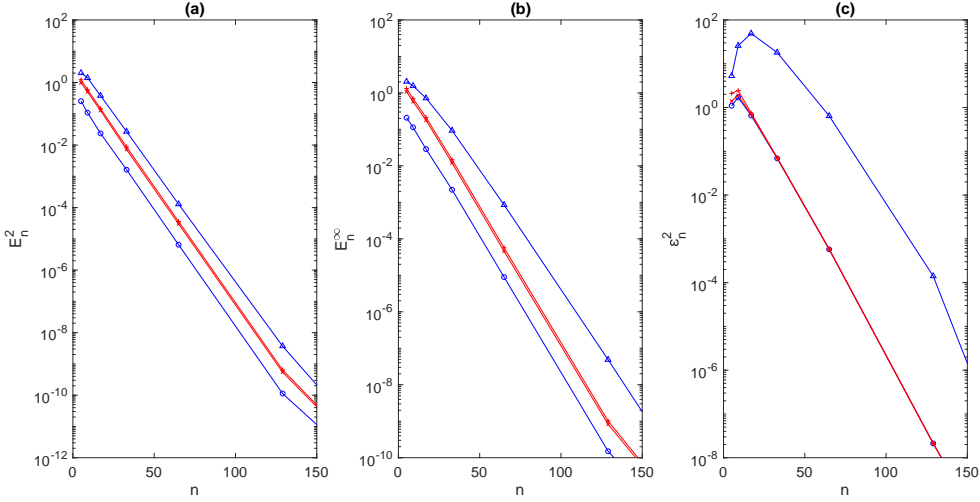}}
\caption{Error estimates (a) $E^{2}_{n}$, (b) $E^{\infty}_{n}$ and (c) $\mathcal{E}^{2}_{n}$ for the sectionally analytic functions $U_1, U_2, L_1, L_2$ (Hurd problem) using the 4-to-1 mapping, as functions of $n$. Parameters are $\chi=\pi/4$, $k=1$, $\theta_0=\pi/3$, $\theta_1=\pi/4$, $\theta_2=\pi/5$. Symbols: $U_1$ ({\color{blue} $\medtriangleup$}), $U_2$ ({\color{blue} $\medcircle$}), $L_1$ ({\color{red} $+$}), $L_2$ ({\color{red} $\times$}).}
\label{Hurdresolution}
\end{figure}

Figure \ref{Hurdresolution} shows the error estimates $E^{2}_{n}$, $E^{\infty}_{n}$, $\mathcal{E}^{2}_{n}$ for $U_1, U_2, L_1, L_2$ using the 4-to-1 mapping, as functions of $n$. We observe spectral convergence, since, as already mentioned, the 4-to-1 mapping is able to capture exactly the far-field behaviour of the sectionally analytic functions. The exact solutions used to compute the error estimates are provided from the highest resolution ($n = 257$) numerical solutions. Although closed-form expressions are given by Hurd \cite{Hurd:1976}, we have not used those to verify our results since their computation is awkward.

\section{Far-field directivity}

The full range Fourier transform was written in the form
\begin{equation}\label{fullrangeFourier}
\Phi(\alpha,y)= \begin{cases}
A(\alpha) \mathrm{e}^{-\gamma y}, \quad & y \geq 0, \\
B(\alpha) \mathrm{e}^{\gamma y}, \quad & y \leq 0.
\end{cases}
\end{equation}
(For the Sommerfeld problem: $B(\alpha)=-A(\alpha)$.) Taking the Fourier inverse of \eqref{fullrangeFourier}, the diffracted field $\phi$ can be written as 
\begin{equation}\label{diffractedfield}
\phi(x,y) = \dfrac{1}{2\pi} {\displaystyle \int_{-\infty}^{\infty}}
{A(\alpha) \ee^{-\ci\alpha x - \gamma y} \dd\alpha}
\end{equation}
in the upper half-plane. To obtain the far-field asymptotics ($x,y \rightarrow \infty$), we deform the integration contour on the steepest descent path \cite{Noble:1988}. Note that the stationary phase method can alternatively be used \cite{Crighton2012modern}.  Application of the steepest descent method to \eqref{diffractedfield} gives
\begin{equation}\label{farasymptotics}
\phi(r,\theta) \sim \dfrac{\sqrt{k} \mathrm{e}^{-{\rm i}\pi/4}}{\sqrt{2\pi}} A(-k \cos \theta) \sin \theta \dfrac{\mathrm{e}^{{\rm i}kr}}{\sqrt{r}}, \quad \text{for} \quad 0 \leq \theta \leq \pi, \quad \text{as} \quad r\rightarrow \infty
\end{equation}
where $(r,\theta)$ are polar coordinates. A similar representation for $\phi(r,\theta)$ in the lower half-plane ($-\pi\leq\theta\leq0$) can be found. The far-field directivity, $D(\theta)$, is defined by \eqref{directivity}. 


Note that the formula \eqref{farasymptotics} breaks down near the shadow boundaries. The reason is that the pole $\alpha=k \cos\theta_0$ is near the saddle point $\alpha=-k\cos\theta$ and, therefore, the steepest descent method which was employed to obtain the far-field asymptotics \eqref{farasymptotics} fails. To overcome this, we can proceed as in \cite{Noble:1988}. Suppose that $A(\alpha)$ can be written as
\begin{equation}\label{Anewexpansion}
A(\alpha)=\frac{H(\alpha)}{\alpha-k\cos\theta_0},
\end{equation}
where function $H(\alpha)$ varies slowly near the saddle point $\alpha=-k \cos\theta$. The form of $H(\alpha)$ depends on the problem under consideration. Then one can break up the solution $\phi(x,y)$ into two terms, one of which is uniformly valid while the other takes the form of a plane wave in the far field. We do not present results from this procedure here, since our goal is to show that the numerical method reproduces the Fourier transform $A(\alpha)$ accurately and we compare to exact results that are not uniformly valid. If desired, one could then use the modification presented in \cite{Noble:1988} to obtain uniformly valid solutions numerically.

\vspace{0.3cm}

\noindent
{\bf The Sommerfeld problem} 
Using \eqref{farasymptotics} and $B(\alpha)=-A(\alpha)$, it follows that the directivity is given by
\begin{equation}\label{directivitynum}
D(\theta) = \frac{\sqrt{k} \mathrm{e}^{-{\rm i}\pi/4}}{\sqrt{2\pi}} A(-k \cos \theta) \sin \theta, \quad \text{for} \quad 0 \leq \theta \leq 2\pi,
\end{equation}
where
\begin{equation}\label{ASom}
A(\alpha)=-\frac{1}{\gamma} \left(\Phi_+'(\alpha,0) + \Phi_-'(\alpha,0)\right)=D_-.
\end{equation}
Expression \eqref{ASom} follows from \eqref{RHSom}. If the argument of $A(\cdot)$ in \eqref{directivitynum} is in the upper half-plane, then we use the first equality in \eqref{ASom} to compute $A$. In this case, the function $\Phi_-'(\alpha,0)$ is defined in \eqref{Phidminus0}, while we use the numerical solution to the Sommerfeld problem above to compute $\Phi_+'(\alpha,0)$. If the argument of $A(\cdot)$ in \eqref{directivitynum} is in the lower half-plane, then we use the latter equality in \eqref{ASom} to compute $A$; in this case, $D_-$ is computed from our numerical solution.

The exact far-field directivity for the hard-hard Sommerfeld problem \cite{Noble:1988} which is used for comparison to our numerical solutions is given by
\begin{equation}\label{DirectivitySom}
D(\theta)=-\sqrt{\frac{2}{k\pi}} \mathrm{e}^{-{\rm i}\pi/4} ~ \frac{\sin(\theta/2) \sin(\theta_0/2)}{\cos \theta+\cos\theta_0}.
\end{equation}

\vspace{0.3cm}

\noindent
{\bf Senior's problem}
Using \eqref{Senior4eqns}, the functions $A(\alpha)$ and $B(\alpha)$ can be expressed in terms of the sectionally analytic functions as
\begin{equation}
A(\alpha)=-\frac{1}{\gamma} \left( \Phi'_+(\alpha,0^+) + \Phi'_{-}(\alpha,0) \right), \qquad B(\alpha)= \frac{1}{\gamma} \left( \Phi'_+(\alpha,0^-) + \Phi'_{-}(\alpha,0) \right)
\end{equation}
and the directivity can be computed using \eqref{farasymptotics}. Similarly to the Sommerfeld problem, different but equivalent representations for the functions $A$ and $B$ using \eqref{systemWH1}--\eqref{systemWH2} can be written which are then used to compute the far-field directivity numerically.

\vspace{0.3cm}

\noindent
{\bf Hurd's problem}
Following Hurd \& Przezdziecki \cite{HurdPrze:1981}, the functions $A(\alpha)$ and $B(\alpha)$ can be expressed as
\begin{equation}
A(\alpha)=\frac{1}{2} \left(L_1+\frac{L_2}{\beta} \right), \qquad B(\alpha)=\frac{1}{2} \left(-L_1+\frac{L_2}{\beta} \right)
\end{equation}
with additional equivalent representations provided by \eqref{Hurdequations1}--\eqref{Hurdequations2}, and directivity can be computed from application of the steepest descent method to the analogous expression \eqref{diffractedfield}, since the definition of the Fourier transform pair in \cite{HurdPrze:1981} is different, with $\gamma$ replaced by $\beta$ (consistent with \eqref{HurdPhi}) and different normalising factor. An analogous representation to \eqref{farasymptotics} can be written with $A(-k \cos \theta)$ replaced by $A(k \cos \theta)$.

Bowman \cite{Bowman:1967} gave the far-field asymptotics for different impedance boundary conditions along the top and bottom sides of the semi-infinite plane: 
\begin{equation}\label{Bowmanfar}
\phi(r,\theta) \sim \frac{1}{4 {\rm i}} \sqrt{\frac{2}{k \pi}} \mathrm{e}^{-{\rm i}\pi/4} U(\theta,\theta_0) \frac{\mathrm{e}^{{\rm i}kr}}{\sqrt{r}}, \quad \text{as} \quad r\rightarrow \infty,
\end{equation}
where $U(\theta,\theta_0)$ is given by
\begin{equation}
U(\theta,\theta_0)=\frac{\sin(\theta_0/2)}{\psi(\pi-\theta_0)} \left[ \frac{\psi(-\theta)}{\sin(\theta/2)+\cos(\theta_0/2)}+\frac{\psi(2\pi-\theta)}{\sin(\theta/2)-\cos(\theta_0/2)} \right],
\end{equation}
where $\psi(x)$ is a known function \cite{Bowman:1967}. Other expressions are given in \cite{HurdPrze:1981,Babich2008diffraction,OsipovNorris:1999}. The expression \eqref{Bowmanfar} is used for comparison to our numerical solutions for Senior's and Hurd's problems.

\begin{figure}
\centering
{\hspace{-3.2em}
\includegraphics[scale=0.3]{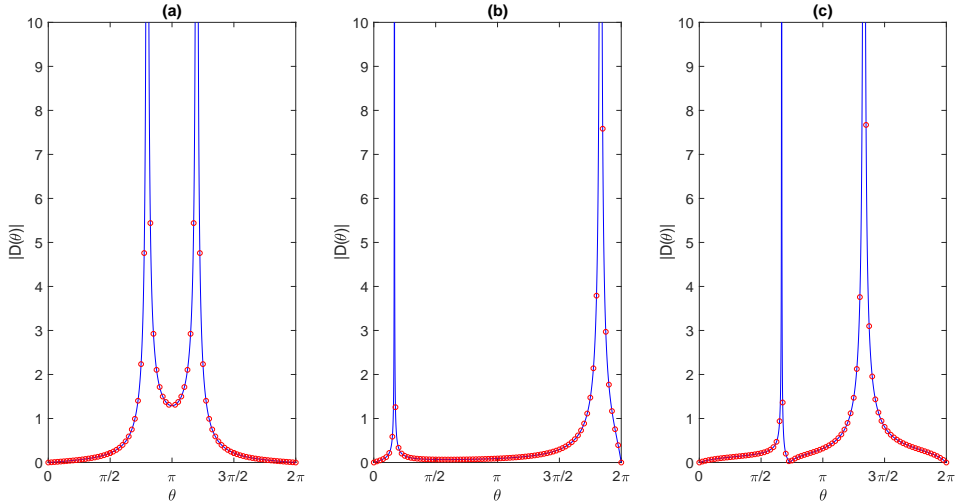}}
\caption{Far-field directivity $|D(\theta)|$ as a function of $\theta$ using exact solutions (solid curves) and numerical solutions using $n = 129$ points (circles). (a) Hard-hard Sommerfeld problem (for $k=1$, $\theta_0=\pi/5$), (b) Senior's problem ($k=1$, $\theta_0=5\pi/6$, $S=\csc(\pi/5)$) and (c) Hurd's problem (for $k=1$, $\theta_0=\pi/3$, $\theta_1=\pi/4$, $\theta_2=\pi/5$).}
\label{fig:directivity}
\end{figure}

Figure \ref{fig:directivity} shows the the far-field directivity $|D(\theta)|$ for the hard-hard Sommerfeld, Senior and Hurd problems, as a function of angle $\theta$, computed using the exact and numerical solutions. The numerical results are indistinguishable from the exact directivities; this reflects the excellent approximation properties of our numerical schemes. We are computing approximations to the Fourier transforms at the stationary phase points. Since the transforms are analytic in their respective half-planes, we obtain spectral convergence. We observe that our numerical method recovers the singularities in the far-field directivities associated with the shadow boundary regions.

\section{Discussion}

We have presented fast and accurate numerical schemes for the solution of scalar and matrix WH problems by exploiting their links with RH problems. The idea is to solve the corresponding RH problems adapting the methods of Olver and Trogdon \cite{Olver:2012,TO:2015} to take into account the known far-field behaviour of the solutions to construct tailor-made numerical schemes. In particular, we have used rational mappings with multiple inverses, extending the M\"obius mappings discussed in \cite{TO:2015} to account for the $\mathcal{O}(\alpha^{-1/2})$ decay of the solutions for large $|\alpha|$. 

The 4-to-1 mapping was used to capture the $\alpha^{-1/2}$ behaviour of the Fourier transforms that stems from the $r^{1/2}$ behaviour near the origin found in diffraction problems with e.g.~hard-hard and impedance-impedance boundary conditions. For soft-hard boundary conditions, $\phi_{j}(r,\theta) \sim r^{1/4} h_{j}(\theta)$, $j=1,2$ and, therefore, the numerical scheme presented here must be adapted (e.g.~using an 8-to-1 rational mapping).

We analysed problems of diffraction of plane waves by a semi-infinite plane with different boundary conditions that produce scalar and matrix WH problems and obtained accurate and spectrally convergent solutions using our numerical schemes. Our results were verified against exact solutions for the cases for which these exist, as well as high-resolution numerical solutions. We  computed the far-field directivity obtained from the far-field asymptotics, which is the most useful physical prediction from these problems.

The numerical scheme presented in this study can be adapted to solve other more complicated problems some of which we mention below. Scattering by wedges and other simple geometric shapes can lead to WH problems (for a discussion of many scattering problems, see \cite{Bowman:1969}), as can linear elasticity problems with straight-line geometry.

Wickham \cite{Wickham:1995} and Abrahams \cite{Abrahams:1997} analysed the problem of scattering of acoustic waves by a semi-infinite screen at the interface between two compressible media with different physical properties, so that $c_{1}\neq c_{2}$, $\rho_{1}\neq \rho_{2}$. The problem was reduced to a matrix WH problem of the form
\begin{equation}\label{AbrahamsmatrixWH}
\left[ \begin{array}{cc} 1 & \mu \delta \\ -\mu/\gamma & 1 \end{array} \right] \Phi_{+}=\Phi_{-}+\frac{2{\rm i} \mu}{\alpha+k \cos\theta_0}  \left[ \begin{array}{c} 0 \\ 1 \end{array} \right],
\end{equation}
where
\begin{equation}
\Phi_{+}=\left[ \begin{array}{c} U_{1} \\ U_{2} \end{array} \right], \qquad \Phi_{-}=\left[ \begin{array}{c} L_{1} \\ L_{2} \end{array} \right]
\end{equation}
and $\gamma(\alpha)=(\alpha^{2}-k^{2})^{1/2}$, $k>1$, $\delta(\alpha)=(\alpha^{2}-1)^{1/2}$, with the same branch cut structure as $\gamma(\alpha)$, and $\mu=\sqrt{\rho_2/\rho_1}$. Local analysis near the origin \cite{Wickham:1995, Abrahams:1997} gives $\phi_{j}(r,\theta) \sim r^{a} h_{j}(\theta)$, $j=1,2$, with $a=\pi^{-1} \tan^{-1}{\mu}$ . Therefore, to construct accurate numerical schemes for solving problems with solutions exhibiting irrational $a$, we need to consider a different approach, e.g.~expanding in Jacobi polynomials rather than the Chebyshev polynomials considered here. This is beyond the scope of this work, but remains a topic for future study. 

The present numerical approach can also be adapted to devise
 accurate numerical schemes for WH problems with jump matrices containing exponential factors \cite{AbrahamsWickham:1990,Antipov:2015,Kisil:2018}, e.g.~of the form
\begin{equation}
\Phi_-=\left[ \begin{array}{cc} A(\alpha) & B(\alpha) \mathrm{e}^{{\rm i}\alpha L} \\ C(\alpha) \mathrm{e}^{-{\rm i} \alpha L} & D(\alpha) \end{array} \right] \Phi_+ + \left[ \begin{array}{c} f_{1}(\alpha) \\ f_{2}(\alpha) \end{array} \right].
\end{equation}
Problems with exponentials $\mathrm{e}^{\pm {\rm i}\alpha L}$ in the elements of the jump matrix arise in acoustics when there are different boundary conditions in the physical space $(-\infty,0)$, $(0,L)$ and $(L,\infty)$, for example in the scattering of acoustic waves by multiple slits or poroelastic plate extensions.

Finally, motivated by various applications such as calculating effective properties of periodic composites \cite{CrasterObnosov:2006}, we aim to adapt the method for solving RH problems with jump matrices defined piecewise along $\mathbb{R}$. To solve such problems, we have to construct numerical schemes and rational mappings with different behaviours at critical points along the contour.

\vspace{0.5cm}

\noindent
{\bf Data Accessibility.} This article has no additional data.

\vspace{0.2cm}

\noindent
{\bf Authors' Contributions.} SGLS had the original idea for this work. The calculations were carried out by both authors equally. Both authors contributed to the writing of the manuscript.

\vspace{0.2cm}

\noindent
{\bf Competing Interests.} The authors declare that they have no competing interests.

\vspace{0.2cm}

\noindent
{\bf Funding.} This research was partly funded by NSF award DMS-1522675.

\vspace{0.2cm}

\noindent
{\bf Acknowledgments.} Conversations with David Abrahams, Sheehan Olver and Tom Trogdon were very helpful.

\bibliography{RHscatter} 
\bibliographystyle{RS}

\end{document}